# THE STADIUM THEOREM

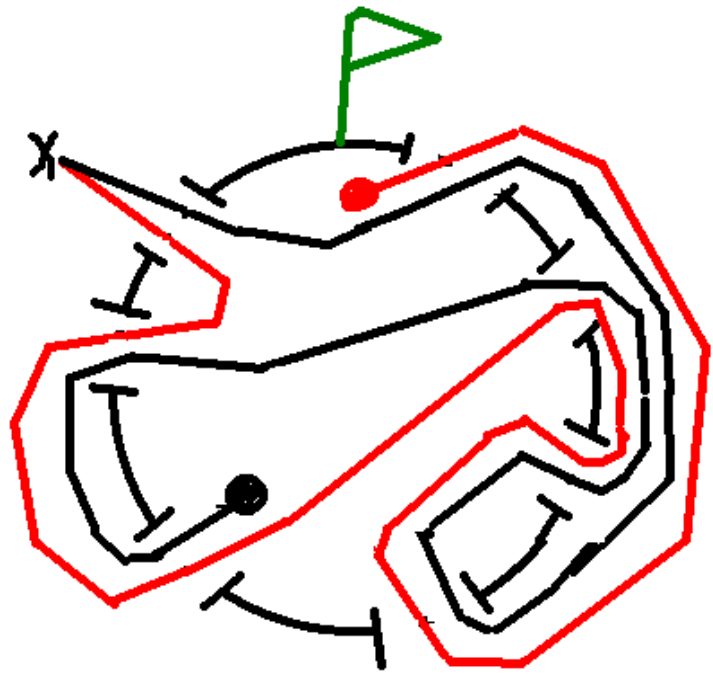
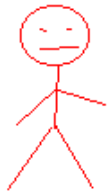
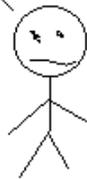

# The Stadium Theorem

This paper gives an informal proof of the **Stadium Theorem** and uses it to prove the odd case of the **General Stadium Theorem**. *The Stadium Theorem.* - Bob and Alice are standing at the same location inside a sports stadium with n gates, n an odd number. Bob says, "I am going to leave this stadium through one gate, come back in through another, go out through a third, etc., and continue in this manner until I pass through all the gates exactly once and I am going to do this in such a way that I do not cross my own path". Alice says, " I'm going to do the same thing and not only not cross my own path but also not cross your path", whereupon Bob says, " then the first gate you go through will have to be the same as the first gate I go through because there are an odd number of gates." The assertion that an odd number of gates implies that the first gate Bob and Alice traverse must be the same is the "Stadium Theorem".

The *General Stadium Theorem* asserts that under the same hypothesis as the *Stadium Theorem* **without assuming an odd number of gates** but with the additional requirement that Bob and Alice arrive at the same location when the walks are completed that **all gates must be traversed by Bob and Alice in <u>exactly</u> the same order.**

The odd case[1] of the *General Stadium Theorem* follows almost immediately from the *Stadium Theorem* as follows. By applying the *Stadium Theorem* to Bob and Alice's paths in both the forward and reverse directions[2] we see that both the first gate and the last gate Bob and Alice traverse must be the same. This fact immediately establishes the case n = 3. A reduction from case n to n-2 for odd n is also easily obtained. A proof of the *General Stadium Theorem* (for n even as well as odd) will be given in a subsequent note.

First we need a few definitions. A *topological n-gon* or simply an *n-gon* is an ordered pair {P ,V}, where P is a simple (i.e., not self-intersecting) closed curve[3] in $R^2$ and $V \subseteq P$, $V = \{ v_1, ... , v_n \}$, is a set of n points on P called the *vertices of P*. Then P-V consists of n mutually disjoint connected pieces of P called the *sides* of P. We also need the facts (Jordan curve theorem) that: $R^2$-P consists of two connected pieces, one bounded piece called the *interior* of P, denoted by *int* P, and one unbounded, called the *exterior* of P, denoted by *ext* P and that any curve with one endpoint in the interior of P and one endpoint in the exterior of P must intersect P.

---

[1] That is, the case extant when the number of gates is odd.
[2] Since Bob and Alice end their walks at the same point, we can consider this point the starting point of two walks in the opposite direction.
[3] All curves we talk about can be taken to be piecewise linear

**The Stadium Theorem.** Let P be an n-gon, n odd, and x be a point, $x \in int$ P and let C = C(x,a), D = D(x,b) denote two simple curves with endpoints x, a and x, b respectively, a ≠ b. That is, C and D are curves which have x as a common endpoint but are otherwise disjoint. Suppose that C(x,a) and D(x,b) each cross each side of P once. Let $Pi_1,...,Pi_n$ and $Pj_1,...,Pj_n$ be the order in which the sides of P are encountered when traversing C(x,a), respectively D(x,b), starting at x. Then we must have $i_1 = j_1$, i.e., each curve crosses the same side of P first.

An informal proof of the theorem is as follows. It suffices to show that given P, C(x,a), D(x,b) which satisfy the hypothesis of the theorem, we can draw a new m-gon P' such that $3 \leq m < n$, where n-m is an even number and such that P', C(x, a), D(x,b) satisfy the hypothesis of the theorem, i.e., such that C(x,a), D(x,b) cross each side of P' once, $C \cap D = x$ and such that the "first" sides that C and D cross are unchanged so that $P'i_1 \neq P'j_1$ if $Pi_1 \neq Pj_1$. Such a redrawing is called a *reduction*. We keep applying reductions until we get the case n=3. We then verify the case n = 3 directly to establish the result.

Let $x_1, ... , x_n$ and $y_1, ... ,y_n$ each be the sequence of points of P encountered when traversing C, respectively D, starting at x. Then each side of P contains exactly one $x_k$ and one $y_j$. Also each *internal section* of C, $C(x_2,x_3), C(x_4,x_5), ... ,C(x_{n-1},x_n)$ and each internal section of D, $D(y_2,y_3), D(y_4,y_5), ... , D(y_{n-1},y_n)$ is contained in the interior of P (except for the endpoints which are on P). Let $S(p,q) := C(x_{2j},x_{2j+1})$ be a given internal section of C. Then P- S consists of two connected pieces $A_1, A_2$ and $A_1 \cup S, A_2 \cup S$ each bound a connected component of  int P - S,  say, $T_1$ and $T_2$. Color the points $x_i$ red and the points $y_k$ black. There are, a priori, three possibilities for an internal section of C as indicated in Figures 1,2, and 3.

The case shown in Figure 1 is ruled out since then $A_1$ would have an odd number of red points while $A_2$ would have an odd number of black points, the former implying $x \in T_1$ while the latter implying $x \in T_2$, a contradiction. The case depicted in Figure 2 gives a desired reduction, yielding a new m-gon with m = n-(e+2) sides (note that $f \geq 3$, otherwise $Pi_1 = Pj_1$ since $x \in T_1$ since $A_1$ has an odd number of reds (and blacks)). The case shown in Figure 3 yields a reduction if e > 0 so let us assume e = 0 for all internal sections of C. The same argument applies to internal sections of D so we can assume e = 0 for all internal sections of D, as well. This is depicted in Figure 4. If we then redraw P as indicated we get a reduction to the case n = 3. (We cannot reduce the case n = 3 any further by this method, since $Pi_1 \neq Pj_1$ would be combined into a single side of P'.) Thus to complete the proof we need to show that the case n = 3 is not possible. But then, referring to Figure 4, s must be connected by an external section to t or t' and u must be connected by an external section to v or v', without the two external sections crossing. That this is not possible may be "intuitively obvious" but to give a rigorous proof of this fact we would need to invoke a deeper theorem of topology. A proof using Kuratowski's theorem, for example, could be constructed by placing a new vertex inside P and connecting it by four simple curves to s, t or t', u, and v or v', appropriately, to get an embedding of $K_5$ in $R^2$, contradicting (the 'only if' part of) Kuratowski's theorem [1]. (Kuratowski's theorem states that a finite graph is planar if and only if it contains no subgraph that is isomorphic to or is a subdivision of $K_5$ or $K_{3,3}$, where $K_5$ is the complete graph of order 5 and $K_{3,3}$ is the complete bipartite graph with 3 vertices in each of the half of the bipartition.)

# References

1. Kazimierz Kuratowski. Sur le problème des courbes gauches en topologie *Fund. Math.*, 15:271-283, 1930.

Figures

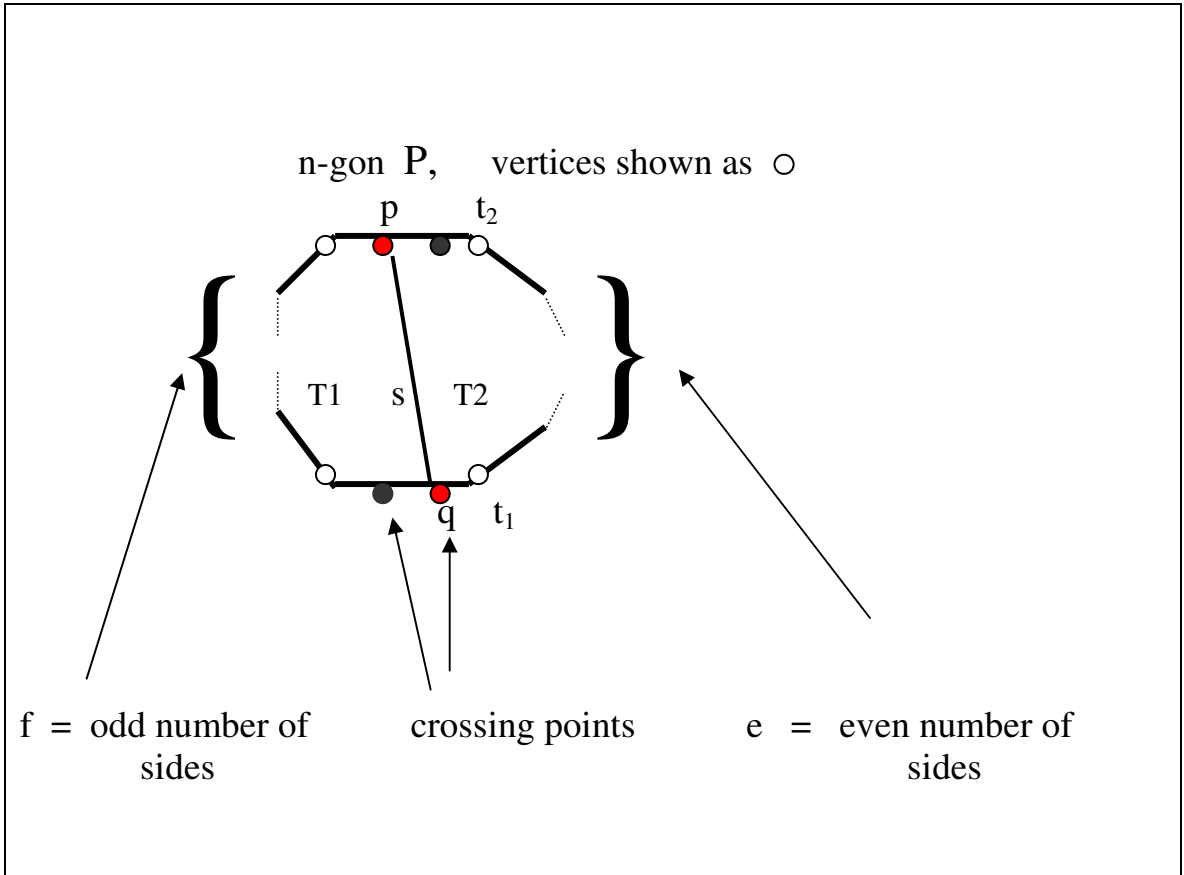

FIGURE 1

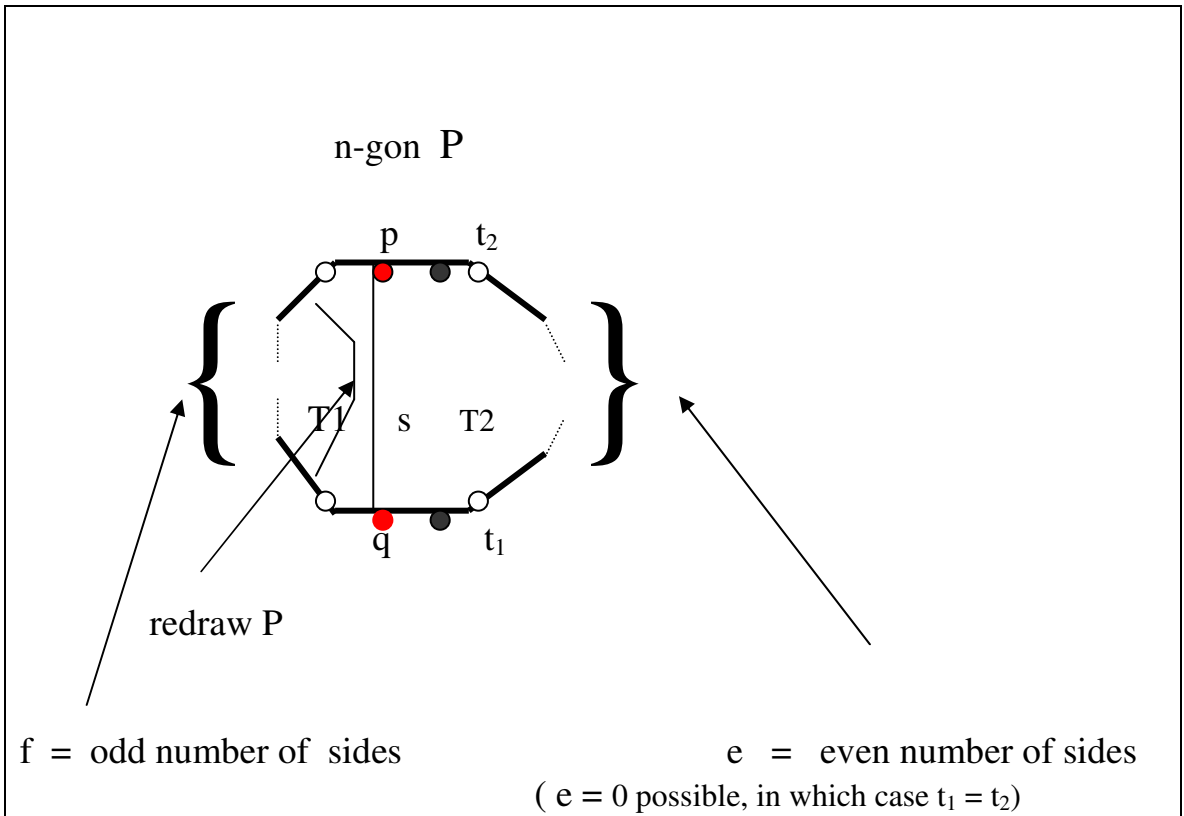

FIGURE 2

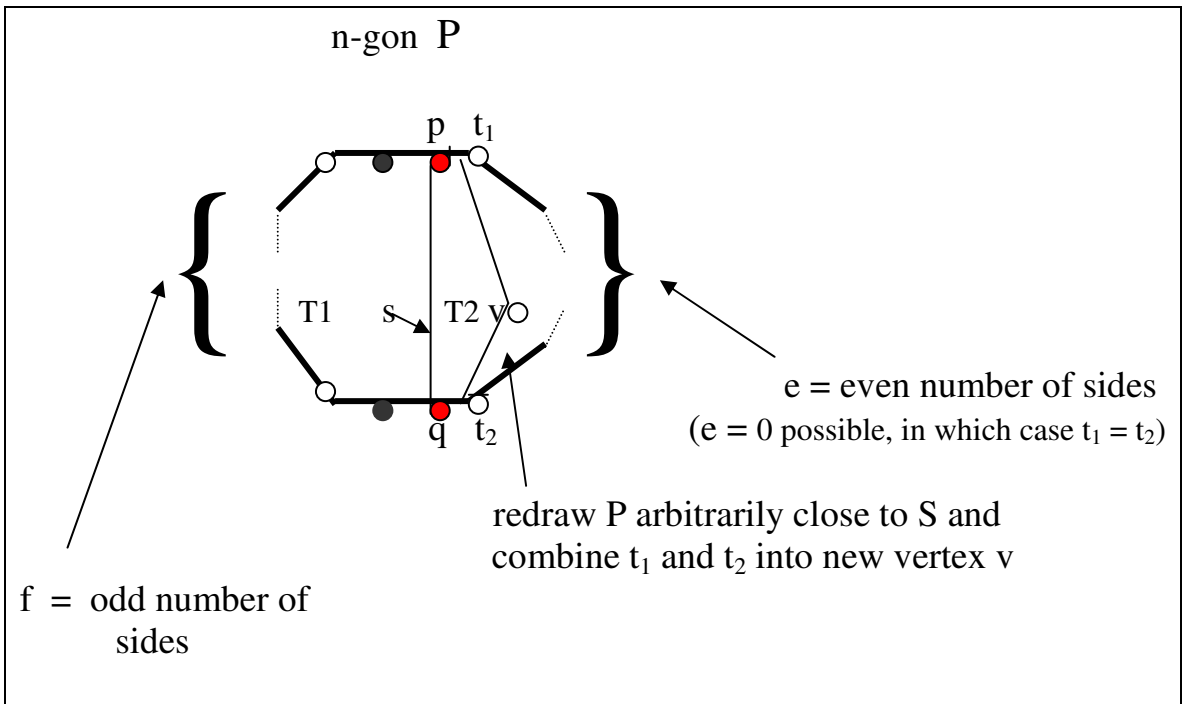

FIGURE 3

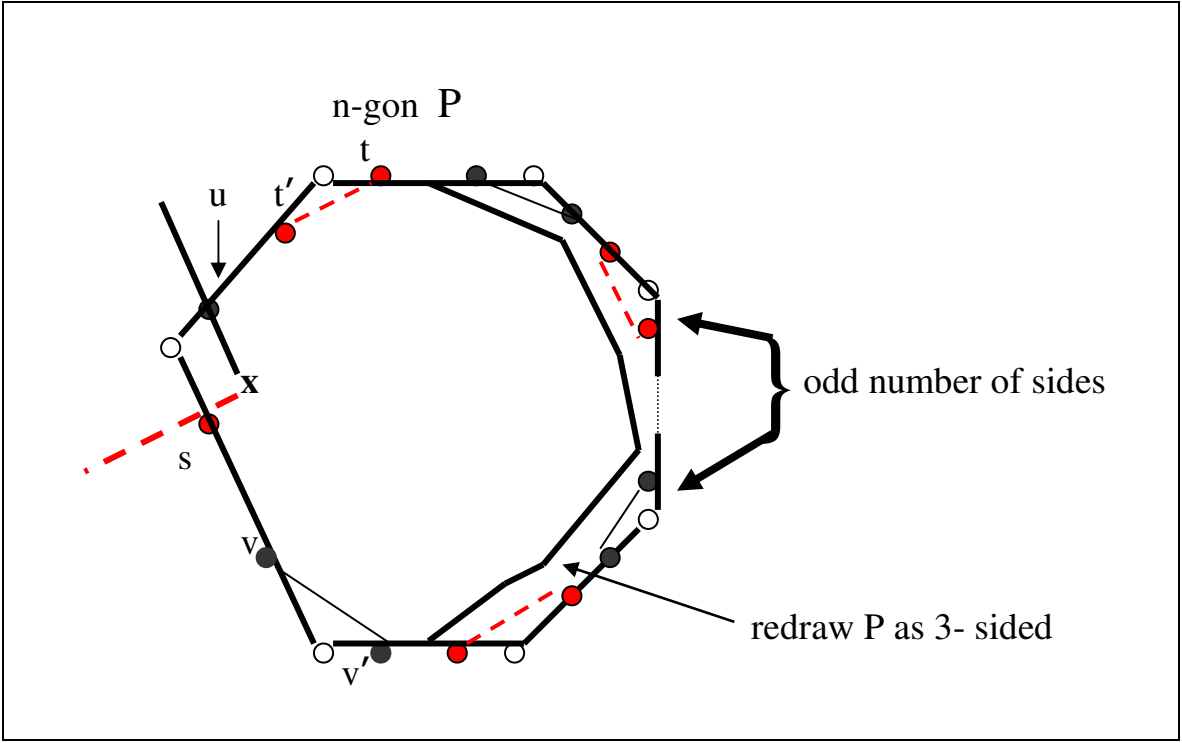

FIGURE 4